\documentstyle[mssymb]{article}
\newcommand{\be}{\begin{equation}}
\newcommand{\ee}{\end{equation}}
\newcommand{\ra}{\rightarrow}
\newcommand{\ul}{\underline}

\newcommand{\geqs}{\geqslant}
\newcommand{\leqs}{\leqslant}

\newcommand{\R}{{\bf R}}
\newcommand{\Z}{{\bf Z}}
\newcommand{\z}{{\bf z}}
\def\wt{\mathop{\rm wt}\nolimits}
\font\svnbf=cmbx7
\newcommand{\sZ}{{\hbox{\svnbf Z}}} 
\newcommand{\sz}{{\hbox{\svnbf z}}} 
\font\fivbf=cmbx5
\newcommand{\tZ}{{\hbox{\fivbf Z}}} 
\begin{document}
\centerline{\Large Lattices and codes with long shadows}
\vspace*{3ex}
\centerline{Noam D. Elkies}
\vspace*{6ex}

\noindent{\bf Introduction.}
\vspace*{1ex}

By a \ul{characteristic} \ul{vector} of an integral unimodular
lattice $L\subset\R^n$ we mean a vector $w\in L$ such that
$(v,w)\equiv(v,v)\bmod2$ for all $v\in L$.  Such vectors
are known to constitute a coset of~$2L$ in~$L$ whose norms
are congruent to~$n$ mod~8 (see e.g.~\cite[Ch.V]{Serre});
dividing this coset by~2 yields a translate of~$L$ called
the \ul{shadow} of~$L$ in~\cite{IEEE}.
If $L=\Z^n$ then $w\in L$ is characteristic if and only if all its
coordinates are odd, so every characteristic vector of~$\Z^n$
has norm at least~$n$.  In~\cite{odd1} we proved that if $L\not\cong\Z^n$
then $L$ has characteristic vectors of norm $\leqs n-8$, and
described without proof all lattices for which $n-8$ is the minimum.
Here we prove this result, and along the way also obtain congruences
and a lower bound on the kissing number of unimodular lattices with
minimal norm~2.  We then state and prove analogues of these results
for self-dual codes, and relate them directly to the lattice
problems via Construction~A.

\vspace*{3ex}

\noindent{\bf Estimates for unimodular lattices}
\vspace*{1ex}

Any integral lattice $L$ decomposes as the direct sum $\Z^r\oplus L_0$
where the $\Z^r$ is generated by the vectors of norm~1 and $L_0$ is
a lattice of minimal norm $\geqs 2$.  [This $L_0$ is called the
``reduced form'' or ``initial lattice'' of~$L$ in~\cite[p.414]{SPLAG},
the latter terminology suggesting the infinite family of lattices
$L_0$, $L_0\oplus\Z$, $L_0\oplus\Z^2$, etc., of which $L_0$ is
the initial member.]  If $L$ (and thus also $L_0$)
is unimodular then the shadow of~$L$ is the
orthogonal sum of the shadows of~$\Z^r$ and~$L_0$.
Replacing $L$ by~$L_0$ thus reduces both the rank of the lattice
and the norm of its shortest characteristic vector by~$r$, and
does not change the difference between these two integers.
We may thus restrict attention to lattices with no vectors of
norm~1 for which that difference is~8, and at the end recover
all such lattices by adding arbitrarily many $\Z$'s.

\vspace*{1ex}

{\bf Theorem 1.}  {\sl
Let $L$ be an integral unimodular lattice in $\R^n$ with no
vectors of norm~$1$.  Then:

i) $L$ has at least $2n(23-n)$ vectors of norm~$2$.

ii) Equality holds if and only if $L$ has no characteristic
vectors of norm $<n-8$.

iii) In that case the number of characteristic vectors
of norm exactly $n-8$ is $2^{n-11} n$.
}

{\sl Proof}\/:  We use theta series as in~\cite{odd1}, though here
we freely invoke modular forms.  For $t$ in the upper half-plane~$H$\/
define
\be
\theta_L(t) := \sum_{v\in L} e^{\pi i |v|^2 t}
= \sum_{k=0}^\infty N_k e^{\pi i k t},
\label{thetadef}
\ee
where $N_k$ is the number of lattice vectors of norm~$k$, and
\be
\theta'_L(t) := \sum_{v\in L+\frac{w}{2}} e^{\pi i |v|^2 t}
= \sum_{k=0}^\infty N'_{\!k} e^{\pi i k t/4},
\label{theta'def}
\ee
where $w\in L$ is any characteristic vector and $N'_{\!k}$ is
the number of characteristic vectors of norm~$k$, or equivalently
the number of shadow vectors of norm~$k/4$.  In~\cite{odd1}
we noted the identity
\be
\theta_L(\frac{-1}{t}+1) = (t/i)^{n/2}\theta'_L(t).
\label{thetaTS}
\ee
By a theorem of Hecke (see e.g.~\cite[Ch.7,~Thm.7]{SPLAG}),
$\theta_L$ is a modular form of weight $n/2$ and can be written as
a weighted-homogeneous polynomial $P_L(\theta_\sZ,\theta_{E_8})$
in the modular forms
\be
\theta_\sZ(t) := 1 +
2 \left( e^{\pi i t} + e^{4\pi i t} + e^{9\pi i t} + \cdots \right) \;
\label{thetaZ}
\ee
of weight $1/2$ and
\be
\theta_{E_8}(t) = 1 + 240 \sum_{m=1}^\infty
 \frac{m^3 e^{2\pi i m t}}{1 - e^{2\pi i m t}}
= 1 + 240 e^{2\pi i t} + 2160 e^{4\pi i t} + \cdots
\ee
of weight 4.  From (\ref{thetaTS}) it follows that $\theta'_L$
is given by
\be
\theta'_L = P_L(\theta'_\sZ,\theta_{E_8}),
\label{theta'L}
\ee
where
\be
\theta'_\sZ(t) = 2\sum_{m=0}^\infty e^{\pi i (m+\frac12)^2 t} =
2 e^{\pi i t/4} \left(
1 + e^{2\pi i t} + e^{6\pi i t} + e^{12\pi i t} + \cdots \right) \; ,
\label{theta'Z}
\ee
and we used the fact that $\theta'_{E_8}=\theta_{E_8}$ because
$E_8$ is an even lattice.  Since $\theta'_\sZ(t) \sim 2e^{\pi i t/4}$
as $t\ra i\infty$, while $E_4(i\infty)=1$ is nonzero, we see
from~(\ref{theta'L}) that the norm of the shortest characteristic
vectors is simply the exponent of~$X$\/ in the factorization
of~$P_L(X,Y)$.\footnote{We could now recover our theorem
from~\cite{odd1} by observing that this exponent is at most~$n$,
with equality if and only if $\theta_L$ is proportional to
$\theta_\tZ^n$, etc.; but this is really the same proof because
the crucial fact that $\theta_\tZ$ vanishes at one cusp and
nowhere else is also an essential ingredient of Hecke's theorem.}

In our setting $N_0=1$ and $N_1=0$.  We first prove part (ii) of
our theorem.  If $L$ has no characteristic vectors of norm $<n-8$
then $P_L(X,Y)$ is a linear combination of $X^n$ and $X^{n-8} Y$.
The known values of $N_0, N_1$ determine this combination uniquely:
we find that
\be
\theta_L = \theta_\sZ^n -
\frac{n}{8} \theta_\sZ^{n-8} \left(\theta_\sZ^8 - \theta_{E_8}\right)
= 1 + 0 \, e^{\pi i t } + 2n(23-n) e^{2\pi i t}
+ \cdots \; .
\label{thetaL}
\ee
Thus $L$ indeed has $2n(23-n)$ vectors of norm~2.  Conversely if
$L$ is an integral unimodular lattice with $N_1=0$ and
$N_2 \leqs 2n(23-n)$ then $n<24$ and $P_L$ has at most 3 terms,
whose coefficients are determined uniquely by $N_0,N_1,N_2$:  
\be
\theta_L = \theta_\sZ^n -
\frac{n}{8} \theta_\sZ^{n-8} \left(\theta_\sZ^8 - \theta_{E_8}\right) +
\frac{N_2 - (2n(23-n))}{16^2}
\theta_\sZ^{n-16} \left(\theta_\sZ^8 - \theta_{E_8}\right)^{\!2}.
\ee
But then by (\ref{theta'L}) we have
\be
N'_{n-16} = 2^{n-24} [N_2 - (2n(23-n))].
\label{n-16}
\ee
Since $N'_{n-16} \geqs 0$ we conclude that $N_2 \geqs 2n(23-n)$
even for $n<24$, as claimed in part~(i) of the theorem; and
equality occurs if and only if $N'_{n-16}$ vanishes, whence
the reverse implication in part~(ii) follows.  Finally to prove
part~(iii) we use (\ref{theta'L},\ref{thetaL}) to compute
\be
\theta'_L = {\theta'_\sZ}^{\!\!n} - \frac{n}{8}
{\theta'_\sZ}^{\!\!n-8} \left({\theta'_\sZ}^{\!\!8} - \theta_{E_8}\right)
= 2^{n-11} n \, e^{(n-8) \pi i t} + \cdots \; ,
\ee
so $N'_{n-8} = 2^{n-11}n$ as claimed.~~$\Box$

Fortunately the integral unimodular lattices of rank $n<24$
are completely known, and those with $N_1=0$ are conveniently listed
with their $N_2$ values in the table of~\cite[pp.416--7]{SPLAG}.
For $n<16$ the shortest characteristic vector must have norm
at least $n-8$, so any unimodular lattice of minimal norm $>1$
must have $2n(23-n)$ vectors of norm~2; this is confirmed
by the table.  When $16 \leqs n \leqs 23$
some lattices can have more than $2n(23-n)$ such vectors, but
it turns out there is always at least one lattice with $N_2=2n(23-n)$.
(Can this be proved a priori?)  Thus, as observed by J.H.~Conway,
the lattices of parts (ii),~(iii) of our theorem are precisely the
integral unimodular lattices of rank $n<24$ with $N_1=0$
that minimize $N_2$ given~$n$.
As noted in~\cite{odd1}, there are fourteen such
lattices; in the following list, adapted from~\cite{odd1},
we label them as in the table of~\cite{SPLAG} by the root system
of {norm-2} vectors:

\vspace*{1ex}

\centerline{
$
\begin{array}{c|cccccccccccc}
n & 8 & 12 & 14 & 15 & 16 & 17 & 18 & 19 & 20 & 21 & 22 & 23
\\ \hline
N_2 & 240 & 264 & 252 & 240 & 224 & 204 & 180 & 152 & 120 & 84 & 44 & 0
\\ \hline
\phantom{\sum^0}\!\!\!\!\! & E_8 & D_{12} & E_7^2 & A_{15} & D_8^2 & A_{11}E_6
 & D_6^3,\ A_9^2 & A_7^2D_5 & D_4^5,\ A_5^4 & A_3^7 & A_1^{22} & O_{23}
\end{array}
$
}

\vspace*{1ex}

We noted in~\cite{odd1} that from our characterization of $\Z^n$
we could also recover the fact that $\Z^n$ is the only integral
unimodular lattice of rank~$n$ for $n<8$.  Likewise from part (iii)
of Theorem~1 we can recover the fact that every integral unimodular
lattice of rank $n<12$ is either $\Z^n$ or $\Z^{n-8}\oplus E_8$.
Indeed there would otherwise be such a lattice of rank~9, 10, or~11
with no vector of norm~1, but then by (iii) the lattice would have
$N'_{n-8}=2^{n-11} n$, which is impossible because $N'_{n-8}$
is an even integer for any lattice of rank $n\neq8$.

Having obtained (\ref{n-16}), we used $N'_{n-16} \geq 0$ to prove
$N_2 \geqs 2n(23-n)$.  Since $N'_{n-16}$ is always an even
integer unless $n=16$ and {\bf0} is a characteristic vector
($\Longleftrightarrow L$ is its own shadow
$\Longleftrightarrow L$ is an even lattice), it follows that in fact
\be
N_2 \equiv 2n(23-n) \bmod 2^{25-n}
\label{25-n}
\ee
for any even unimodular lattice with no vectors of norm~1,
with the exception of the two even lattices $E_8^2$, $D_{16}^+$
of rank~16, which have $N_2=480$.  (This is similar to the
argument used in~\cite[Ch.19,~p.440]{SPLAG} to prove that there are
no ``extremal Type~I lattices'' with $16 \leqs n \leqs 22$,
a fact which now also follows from part~(i) of our Theorem.)
The congruence~(\ref{25-n}) is confirmed by the table of~\cite{SPLAG},
which also reveals that $N_2 - 2n(23-n)$ is always a multiple
of~$2^4$ even for $n=22$ and $n=23$; a ``conceptual'' (but far from
easy) proof of this is found in \cite[Thm.~4.4.2(3)]{reb:thesis}.
Note that even though we have only proved (\ref{25-n}) for $n<24$,
it in fact holds for all~$n$, since $N_2$ is always an even integer.

\vspace*{3ex}

\noindent{\bf Estimates for self-dual binary codes}
\vspace*{1ex}

We recall some basic facts about binary linear codes;
see for instance \cite[Ch.3,~\S2.2]{SPLAG}.
Let $F=\Z/2\Z$ be the two-element field.  We work in the vector
space~$F^n$, whose elements we regard as ``words'' of length~$n$
whose ``letters'' are taken from the ``alphabet''~$F$.  The (Hamming)
``weight'' $\wt(w)$ of a word $w=(w_1,w_2,\ldots,w_n)\in F^n$ is
$\#\{j: w_j = 1\}$ of nonzero coordinates of~$w$.  A ``binary
linear code'' of length~$n$ is a subspace $C \subset F^n$.
A binary \ul{self-dual} code is a linear code which is its own
annihilator under the nondegenerate pairing $(\cdot,\cdot):
C \times C \ra F$, defined by $(v,w)=\sum_{j=1}^n v_j w_j$.
Such a code must have dimension $n/2$, and thus can only exist
if $n$ is even, which we henceforth assume.  Note that under
our pairing we have
\be
(w,w) = (w,1^n) \equiv \wt(w) \bmod 2
\label{mod2}
\ee
for all $w\in F^n$, where $1^n$ is the all-ones vector in~$F^n$.
Thus if $C$\/ is a self-dual code then $C\ni 1^n$ and all the
words in~$C$\/ have even weight.

The ``weight enumerator'' $W_C$ of~$C$\/ is a generating
function for the weight distribution of~$C$\/:
\be
W_C(x,y) := \sum_{c\in C} x^{n-\wt(c)} y^{\wt(c)}.
\label{Wdef}
\ee
For a binary self-dual code a theorem of Gleason (Thm.6 of
\cite[Ch.7]{SPLAG}), analogous to Hecke's theorem for theta
series of lattices, states that $W_C$ is a weighted-homogeneous
polynomial $P_W(x^2+y^2, x^8+14x^4y^4+y^8)$ in the weight enumerators
of the double repetition code $\z := \{(0,0),(1,1)\}\subset F^2$ and
the extended Hamming code in~$F^8$ respectively.

Analogous to the homomorphism $v\mapsto|v|^2\bmod2$ from an integral
lattice to~$\Z/2$ we have for any self-dual code $C\subset F^n$
a linear map from~$C$\/ to~$F$\/ taking any $c\in C$\/ to
$\frac12 \wt(c) \bmod 2$.  We can use the pairing on~$F^n$ to
represent any linear functional on a self-dual code by a unique
coset of the code; thus we find a coset~$C'$ of~$C$\/ consisting
of all $c'\in F^n$ such that
\be
\frac12 \wt(c) \equiv (c,c') \bmod 2
\label{mod4}
\ee
for all $c\in C$.  As in~\cite{IEEE} we call $C'$ the \ul{shadow}
of~$C$, in analogy with the shadow of an integral unimodular lattice.
Let
\be
W'_C(x,y) := \sum_{c\in C'} x^{n-\wt(c)} y^{\wt(c)}
\label{W'def}
\ee
be the generating function for the weight distribution of~$C'$.
Using discrete Poisson inversion as in the proof of the MacWilliams
identity and the characterization~(\ref{mod4}) of~$C'$ we find
as in~\cite{IEEE}
$$
W'_C(x,y) = 2^{-n/2} \sum_{c\in C}
(-1)^{\wt(c)/2} (x+y)^{n-\wt(c)} (x-y)^{\wt(c)}
$$ \vspace*{-4ex} \be \phantom. \label{W&W'} \ee \vspace*{-4ex} $$
= 2^{-n/2} W_C(x+y,i(x-y)).
$$
Thus from $W_C(x,y)=P_W(x^2+y^2, x^8+14x^4y^4+y^8)$ we obtain
\be
W'_C(x,y) = P_W(2xy, x^8+14x^4y^4+y^8) .
\label{W'_C}
\ee
Note that all the words in the shadow thus
have weight congruent to~$n/2$ mod~4.  We could have also obtained
this directly from the MacWilliams identity
\be
W_C(x,y) = 2^{-n/2} W_C(x+y,x-y),
\label{McW}
\ee
(which also underlies Gleason's theorem); this would more closely
parallel the analytic proof of $|w|^2\equiv n\bmod 8$ in~\cite{odd1}.

If $c\in C$\/ has weight~2 then every codeword either contains or
is disjoint from~$c$.  Thus $C$\/ decomposes as a direct sum of
a double repetition code~$\z$ generated by~$c$ and the self-dual
code of length~$n-2$ consisting of codewords disjoint from~$c$.
Iterating this we decompose~$C$\/ as $C_0\oplus\z^r\!$, where
$r$ is the number of weight-2 words in~$C$, and $C_0$ is a self-dual
code of length $n-2r$ with no words of weight~2.  Now for any
self-dual codes $C_1,C_2$, their direct sum $C_1\oplus C_2$
has shadow
\be
(C_1 \oplus C_2)' = C'_1 \oplus C'_2 \, .
\label{oplus}
\ee
Since the shadow of $\z$ is $\{(0,1),(1,0)\}$
it follows that the shadow of~$\z^r$ consists
entirely of words of weight~$r$, and if $C=C_0\oplus\z^r$
then the minimal weight of $C'_0$ is $r$ less than that of~$C'$.

Since $C'$ contains $w+1^n$ whenever it contains $w$ it is clear
that the minimal weight of $C'$ cannot exceed the value $n/2$
attained by $\z^{n/2}$.  This is much easier than proving the
corresponding fact for characteristic vectors of unimodular
lattices, but it does not show that $\z^{n/2}$
is the only self-dual code whose shadow has
minimal weight~$n/2$.  We prove this, as we did for lattices,
by noting that such a code~$C$\/ must have $W'_C(x,y)=(2xy)^{n/2}$,
whence $W_C(x,y)=(x^2+y^2)^{n/2}$.  Since $C$\/ contains $n/2$ words
of weight~2, then, it can only be~$\z^{n/2}$.

We have shown that the shadow of a binary self-dual
code~$C$\/ other than $\z^{n/2}$ contains some words of weight $<n/2$.
Thus $C'$ has minimal weight at most $(n-8)/2$.  We next characterize
all~$C$\/ attaining this bound.  If $C=C_0\oplus\z^r$ then $C$\/
attains the bound if and only if $C_0$ does, so we need only consider
codes without weight-2 words.

\vspace*{1ex}

{\bf Theorem 1A.}  {\sl
Let $C$ be a binary self-dual code of length~$n$ with no
codewords of weight~$2$.  Then:

i) $C$ has at least $n(22-n)/8$ codewords of weight~$4$.

ii) Equality holds if and only if the shadow
of~$C$ contains no codewords of weight $<(n-8)/2$.

iii) In that case the number of codewords of weight exactly
$(n-8)/2$ in the shadow is $2^{(n-14)/2} n$.
}

{\sl Proof}\/: We can mimic the proof of Theorem~1.
If the minimal weight of~$C'$ is at least $(n-8)/2$
then $W'_C$ is a linear combination of $(xy)^{n/2}$
and $(xy)^{(n-8)/2}(x^8+14x^4y^4+y^8)$, and thus $W_C(x,y)$
is a linear combination of $(x^2+y^2)^{n/2}$ and
$(x^2+y^2)^{(n-8)/2}(x^8+14x^4y^4+y^8)$.  The condition
that $C$\/ have no weight-2 codewords then forces
$$
W_C(x,y) = (x^2+y^2)^{n/2} - \frac{n}{8} (x^2+y^2)^{(n-8)/2}
\left( (x^2+y^2)^4 - (x^8+14x^4y^4+y^8) \right) 
$$ \vspace*{-3.5ex} \be \phantom. \label{WC} \ee \vspace*{-6.5ex} $$
= x^n + 0 \, x^{n-2} y^2 + \frac{n(22-n)}{8} x^{n-4} y^4 + \cdots
$$
and
$$
W_C(x,y) = (2xy)^{n/2} - \frac{n}{8} (2xy)^{(n-8)/2}
\left[ (2xy)^4 - (x^8+14x^4y^4+y^8) \right]
$$ \vspace*{-5ex} \be \phantom. \label{W'C} \ee \vspace*{-5ex} $$
= 2^{(n-14)/2} (xy)^{(n-8)/2}
\left[ n x^8 + (128-2n) x^4y^4 + n y^8 \right] .
$$
If $n<24$ and $C$\/ is any binary self-dual code of length~$n$
containing no words of weight~2 and  $n(22-n)/8 + d$\/ words
of weight~4 then its weight enumerator exceeds (\ref{WC}) by
\be
\frac{d}{16} (x^2+y^2)^{(n-16)/2}
\left[x^8+14x^4y^4+y^8 - (x^2+y^2)^4\right]^2,
\label{WC+}
\ee
so the weight enumerator of the shadow~$C'$ exceeds
(\ref{W'C}) by
\be
\frac{d}{16} (2xy)^{(n-16)/2} (x^8-2x^4y^4+y^8)^2.
\label{W'C+}
\ee
Thus $C'$ contains $2^{(n-24)/2} d$\/ words of weight $(n-16)/2$,
from which we find that $d$\/ is a nonnegative multiple of
$2^{(24-n)/2}$.

Alternatively we could deduce Theorem~1A from Theorem~1 via
Construction~A \cite[Ch.7,~\S2]{SPLAG}.  Recall that this construction
associates to a self-dual code $C\subset F^n$ the unimodular integral
lattice
\be
L_C := \{ 2^{-1/2} v \; | \;  v \in \Z^n, v \bmod 2 \in C \}.
\label{LC}
\ee
The theta series of this lattice is given by
\be
\theta_L(t)=W_C(\theta_\sZ(2t),\theta'_\sZ(2t));
\label{thetaLC}
\ee
in particular $L_C$ has no vectors of norm~1 if and only if
$C$\/ has no codewords of weight~2 (NB $L_\sz\cong\Z^2$),
and the $N_2(L_C)-2n$ is $2^4$ times the number of weight-4
codewords in~$C$.  Moreover the set of characteristic
vectors of~$L_C$ is
\be
\{ 2^{1/2} v \; | \;  v \in \Z^n, v \bmod 2 \in C' \}
\label{L'C}
\ee
(in effect the shadow of~$L_C$ is obtained by applying
Construction~A to the shadow of~$C$\/),
so the norm of the shortest characteristic vectors is half
the minimal weight of~$C'$.  Applying Theorem~1 to $L_C$
thus yields Theorem~1A immediately.~~$\Box$

Thus also the codes~$C$\/ of parts (ii),~(iii) of Theorem~1A
are precisely those for which $L_C$ is one of the 14 lattices
listed in connection with Theorem~1.  Of course not every such
lattice arises because $n$ must be even; moreover, the root
system can only involve $A_1$, $D_{2m}$, $E_7$ and $E_8$ if the
lattice arises from construction~$A$.  This leaves only
the seven lattices with root systems
$E_8, D_{12}, E_7^2, D_8^2, D_6^3, D_4^5, A_1^{22}$ of rank
$8, 12, 14, 16, 18, 20, 22$ respectively.  It turns out
that each of those lattices arises as $L_C$ for a unique
code~$C$~\cite{codes1,codes2}.
For instance the first of these arises from the extended Hamming code,
and the last from what might be called the shorter binary Golay code;
these are the shortest self-dual binary codes having minimal
weight~4 and~6 respectively.  Again it so happens that whenever
there is a self-dual code of length~$n<24$ with minimal weight at
least~4, there is such a code (this time unique) with only $n(22-n)/8$
words of weight~4, so Conway's description of our fourteen lattices
also applies {\sl mutatis mutandis} to our seven codes.

\vspace*{3ex}

\noindent{\bf Can we go past $n-8$?}
\vspace*{1ex}

\nopagebreak 

Our results suggest the following questions:

\nopagebreak 

For any $k>0$ is there $N_k$ such that
every integral unimodular lattice all of whose
characteristic vectors have norm $\geq n-8k$ is of the form
$L_0 \oplus \Z^r$ for some lattice $L_0$ of rank at most $N_k$?

\pagebreak 

For any $k>0$ is there $n_k$ such that
every binary self-dual code of whose shadow has
minimal norm $\geq(n-8k)/2$ is of the form
$C_0 \oplus \z^r$ for some code $C_0$ of length at most $n_k$?

\nopagebreak 

Of course a positive answer for lattices would imply one for
codes, and vice versa for a negative answer, via Construction~A,
with $n_k\leq N_k$ in the former case.  Even $k=2$ seems difficult.

\vspace*{1ex}

{\bf Acknowledgements.} 
Thanks to John H. Conway for helpful correspondence.
This work was made possible in part by funding from
the National Science Foundation and the Packard Foundation.

\vspace*{5ex}
\begin{small}
\noindent 
Dept.\ of Mathematics\\
Harvard University\\
Cambridge, MA 02138 USA\\
e-mail: {\tt elkies@zariski.harvard.edu}\\ \\
February, 1995
\end{small}

\end{document}